\documentclass[11pt]{article}

\usepackage{amsmath, amssymb, amsthm}   
\usepackage{mathrsfs}                  
\usepackage{graphicx}                  
\usepackage{enumitem}                  
\usepackage{hyperref}                  
\usepackage{xcolor}                    

\usepackage{lmodern}                  
\usepackage[T1]{fontenc}
\usepackage[utf8]{inputenc}

\newtheorem{theorem}{Theorem}[section]

\title{Algebraic Constructions of Universal Cycles on Grassmannians $G_q(2,n)$}
\author{Chen-Yu Chi, Ming Hsuan Kang,  Yu Hsuan Hsieh, }
\date{June 2025}

\begin{document}

\maketitle

\begin{abstract}
We study universal cycles on the Grassmannian $G_q(2,n)$, the set of $2$-dimensional $\mathbb{F}_q$-subspaces of $\mathbb{F}_q^n$. 
While their existence is known from inductive and Eulerian graph methods, we give a direct algebraic construction when $n$ is odd under the coprimality condition $\gcd(n,\,q(q^2-1))=1$, using a projective-ratio decomposition and a global product condition. 
We also present explicit examples where a single cycle is simultaneously universal for both $G_q(2,5)$ and $G_q(3,5)$, realizing Grassmannian duality $|G_q(k,n)|=|G_q(n-k,n)|$ at the level of universal cycles.
\end{abstract}

\section{Introduction}

The classical De Bruijn sequence provides an elegant cyclic listing of all strings of fixed length over a finite alphabet, with each substring appearing exactly once. First introduced by Flye Sainte-Marie and later popularized by De Bruijn \cite{debruijn1946circuits}, such sequences have had lasting influence in combinatorics, coding theory, pseudorandomness, and universal cycles \cite{fredricksen1982survey,rothemund1996program,ruskey2001survey}.

Recent work has sought to extend this paradigm to algebraic and geometric contexts, replacing strings with vector configurations and substrings with subspaces. In particular, the Grassmannian
\[
G_q(2,n) = \{ \text{2-dimensional $\mathbb{F}_q$-subspaces of $\mathbb{F}_q^n$}\}
\]
serves as a natural generalization of the classical setting. A long-standing problem is to construct a cyclic sequence of vectors
\[
v_0, v_1, \dots \in \mathbb{F}_q^n
\]
such that each consecutive pair $\{v_i,v_{i+1}\}$ spans a distinct element of $G_q(2,n)$, and every such subspace appears exactly once. The existence of such universal cycles for $G_q(2,n)$ is already known, established by recursive and inductive methods combined with Eulerian graph techniques \cite{jackson2009recursive}. However, these methods are not explicit.

It is worth noting that the parity of $n$ plays a fundamental role. Already in the classical case of $2$-subsets of $\{1,\dots,n\}$ (which can be regarded as the degenerate case $q=1$), a universal cycle exists only when $n$ is odd, since an Eulerian circuit in $K_n$ requires even degree at each vertex. Thus the restriction to odd $n$ in our algebraic construction is a natural analogue of this classical phenomenon.

Our goal is to revisit this problem from a purely algebraic perspective and to give a direct construction when $n$ is odd under a natural coprimality condition.  
The method is efficient and simple: it requires only the choice of a primitive polynomial for $\mathbb{F}_{q^n}$ and a set of projective-ratio representatives.  
Moreover, every ordering of these representatives yields a valid universal cycle for $G_q(2,n)$, providing flexibility to impose additional properties.  
In particular, we give explicit small examples where the same cycle is simultaneously universal for both $G_q(2,5)$ and $G_q(3,5)$, showing that Grassmann duality $|G_q(k,n)|=|G_q(n-k,n)|$ can, in special cases, be realized directly at the level of universal cycles.

\section{Algebraic Decomposition of $G_q(2,n)$}
\label{sec:algebraic-decomposition}

\subsection{Projective Ratios}

We identify $\mathbb{F}_q^n$ with the extension field $E = \mathbb{F}_{q^n}$ and write $F = \mathbb{F}_q$. For $W = \mathrm{span}_F\{v,w\}$ a 2-subspace, define the \emph{projective ratio}
\[
\operatorname{pr}(v,w) := v/w \in E^\times \setminus F^\times,
\]
well-defined up to the action of $\mathrm{PGL}_2(F)$ via Möbius transformations. This yields a projection
\[
\Phi : G_q(2,n) \to \mathcal{C} := (E^\times \setminus F^\times)/\mathrm{PGL}_2(F),
\]
whose fibers $\Phi^{-1}([c])$ correspond to subspaces with the same projective ratio class $[c]$.

\subsection{Fibers and Non-Collapsing Condition}

Each fiber has size $|\Gamma| = \tfrac{q^n-1}{q-1}$, where $\Gamma = E^\times/F^\times$.  
Note that $\Gamma$ can be naturally identified with the Grassmannian $G_q(1,n)$, since cosets of $F^\times$ in $E^\times$ correspond exactly to $1$-dimensional $F$-subspaces of $E$.  
When $n$ is odd, $E$ has no quadratic subfield, so every fiber is uniform of this size. Consequently, to enumerate all 2-subspaces it suffices to select a set of representatives with distinct projective ratios, and then generate the entire fiber by multiplying these representatives by $\alpha^i$, where $\alpha$ is a generator of $\Gamma$.

Let $\sigma(z)=z^q$ denote Frobenius, generating $\mathrm{Gal}(E/F) \cong \mathbb{Z}/n\mathbb{Z}$. Since the Galois and Möbius actions commute, we define the \emph{collapse degree}
\[
m_z := \big| \mathrm{PGL}_2(F)\cdot z \ \cap \ \langle \sigma \rangle \cdot z \big|.
\]
We say the action is non-collapsing if $m_z=1$ for all $z$. A sufficient condition is
\[
\gcd(n,\, q(q^2-1))=1,
\]
since two commuting group actions with coprime orders have transversal orbits.

\section{The Global Product Construction}

Under the non-collapsing condition, each $\mathrm{PGL}_2(F)$-orbit intersects each Galois orbit in exactly one point. Partition
\[
\mathcal{C} = C_1 \sqcup \cdots \sqcup C_m, \qquad
C_i = \{ [g_i], [g_i^q], [g_i^{q^2}], \dots\},
\]
for chosen representatives $g_1,\dots,g_m$. More generally, one can choose $g_1 \in E^\times$ such that $\gamma g_1 = \alpha g_1$ for some $\gamma \in \mathrm{PGL}_2(F)$.  
As a concrete example, taking 
\[
g_1 = z = \tfrac{1}{\alpha-1},
\]
we have
\[
\begin{pmatrix} 1 & 1 \\ 0 & 1 \end{pmatrix} z = \alpha z.
\]

We then define the representative system
\[
\{c_1,\dots,c_r\} = \{\alpha g_1, g_1^q,g_1^{q^2},\dots\}
\sqcup \{g_2,g_2^q,g_2^{q^2},\dots\} \sqcup \cdots \sqcup \{g_m,g_m^q,g_m^{q^2},\dots\}.
\]
The product satisfies
\[
c_1 c_2 \cdots c_r \in \alpha F^\times.
\]

\section{Constructing the Cycle}

Since \(n\) is odd, the extension field \(E=\mathbb{F}_{q^n}\) contains no quadratic subfield, and hence every fiber of  
\[
\Phi : G_q(2,n) \to \mathcal{C}
\]  
has uniform size  
\[
|\Gamma| = \frac{q^n-1}{q-1}, \qquad \Gamma = E^\times/F^\times.
\]
Thus every 2-subspace of \(\mathbb{F}_q^n\) can be generated by first choosing a representative with a distinct projective ratio class, and then multiplying it by successive powers of a generator \(\alpha \in \Gamma\).  

Formally, define a sequence \(\{\beta_i\}\subset E^\times\) by  
\[
\beta_0 := 1, 
\qquad 
\beta_i := c_1 c_2 \cdots c_i ,
\]  
where the indices of \(c_i\) are taken modulo \(r\). Then the ratios satisfy  
\[
\frac{\beta_i}{\beta_{i-1}} = c_i, 
\qquad 
\beta_{i+r} = \alpha \beta_i \pmod{F^\times}.
\]

Finally, associate to each \(\beta_i\) the 2-subspace  
\[
W_i := \mathrm{span}_F\{\beta_i, \beta_{i-1}\}.
\]

This construction yields a cyclic sequence \(\{W_i\}\) that traverses all 2-subspaces in \(G_q(2,n)\) exactly once, with periodicity given by multiplication by \(\alpha\).

This procedure produces a cyclic sequence of 2-subspaces in $G_q(2,n)$.  
Two useful structural properties follow:  
\begin{enumerate}
    \item Uniformity: every element of $G_q(1,n)$ (i.e. every line through the origin in $\mathbb{F}_q^n$) appears in the cycle the same number of times, reflecting the uniformity of fibers.  
    \item Permutation invariance: the order of the representatives $c_1,\dots,c_r$ does not affect the universal cycle property.  
   Hence, by permuting these factors one can construct alternative universal cycles, sometimes with additional desirable structure (for example, simultaneously realizing universality in both $G_q(2,n)$ and $G_q(n-2,n)$).
\end{enumerate}

\begin{theorem}[Universal Cycles for Odd $n$]
Let $n$ be odd and $q$ a prime power such that $\gcd(n,\, q(q^2-1))=1$. 
Then the sequence constructed from the cyclic product system
\[
\beta_0 := 1, \qquad 
\beta_i := c_1c_2\cdots c_i, \qquad
W_i := \mathrm{span}_F\{\beta_i,\beta_{i-1}\},
\]
where $\{c_1,\dots,c_r\}$ are chosen as in Section~\ref{sec:algebraic-decomposition} and the indices of \(c_i\) are taken modulo \(r\),
is a universal cycle on $G_q(2,n)$. 
That is, every 2-dimensional subspace of $\mathbb{F}_q^n$ appears exactly once in the sequence, with periodicity
\[
W_{i+r} = \alpha \cdot W_i
\]
for a generator $\alpha \in \Gamma$.
\end{theorem}

\paragraph{Remark (Even $n$).}  
The restriction to odd $n$ is natural. In the classical case of $2$-subsets of $\{1,\dots,n\}$ (the degenerate case $q=1$), a universal cycle exists only when $n$ is odd, since an Eulerian circuit in $K_n$ requires each vertex to have even degree. 
In the finite field setting, the obstruction manifests as subfield planes from $\mathbb{F}_{q^2}\subset \mathbb{F}_{q^n}$ when $n$ is even. 
These planes force certain projective-ratio classes to repeat, and the construction then yields only an “almost” universal cycle in which each subfield plane appears $q+1$ times.

\section{Example: $G_2(2,5)$}

Let $q=2$, $n=5$, so $E=\mathbb{F}_{2^5}$ with primitive element $\alpha$ a root of $x^5+x^2+1$. Then $E^\times$ is cyclic of order $31$, and since $F^\times=\{1\}$ we have $\Gamma = E^\times$.

The Möbius action of $\mathrm{PGL}_2(F)\cong S_3$ partitions $\Gamma$ into 5 orbits of size 6. Explicitly, these are
\[
\begin{aligned}
&\{\alpha^1, \alpha^{13}, \alpha^{14}, \alpha^{17}, \alpha^{18}, \alpha^{30}\}, \\
&\{\alpha^2, \alpha^3, \alpha^5, \alpha^{26}, \alpha^{28}, \alpha^{29}\}, \\
&\{\alpha^4, \alpha^6, \alpha^{10}, \alpha^{21}, \alpha^{25}, \alpha^{27}\}, \\
&\{\alpha^7, \alpha^9, \alpha^{15}, \alpha^{16}, \alpha^{22}, \alpha^{24}\}, \\
&\{\alpha^8, \alpha^{11}, \alpha^{12}, \alpha^{19}, \alpha^{20}, \alpha^{23}\}.
\end{aligned}
\]

Under Frobenius $x \mapsto x^2$, these five Möbius orbits are cyclically permuted, so together they form a single Galois orbit. Thus it suffices to choose one representative $g_1$ for the entire system. Take
\[
g_1 = \alpha^2, \qquad c_1 = g_1 \alpha = \alpha^3.
\]
The remaining representatives are then the Frobenius conjugates of $g_1$:
\[
c_2 = g_1^2 = \alpha^4, \quad
c_3 = g_1^{2^2} = \alpha^8, \quad
c_4 = g_1^{2^3} = \alpha^{16}, \quad
c_5 = g_1^{2^4} = \alpha^{32} = \alpha.
\]

Altogether, the cyclic product system is
\[
\{c_1,c_2,c_3,c_4,c_5\} = \{\alpha^3, \alpha^4, \alpha^8, \alpha^{16}, \alpha\}.
\]
Their product is
\[
c_1 c_2 c_3 c_4 c_5 = \alpha^{3+4+8+16+1} = \alpha^{32} = \alpha,
\]
a generator of $\Gamma$. 

Defining
\[
\beta_0 := 1, \qquad \beta_i := c_1 c_2 \cdots c_i,
\]
where the indices of \(c_i\) are taken modulo \(r\)
and defining
\[
W_i := \mathrm{span}_F\{\beta_i, \beta_{i-1}\},
\]
we obtain a cyclic sequence of $155$ distinct subspaces in $G_2(2,5)$, each appearing exactly once. 

Thus this explicit construction realizes the theorem in the case $q=2$, $n=5$, and exhibits a universal cycle for $G_2(2,5)$.

\paragraph{Additional Remark.}
By Grassmann duality, $|G_q(k,n)| = |G_q(n-k,n)|$, so it is natural to ask whether a single universal cycle may simultaneously serve both $G_q(k,n)$ and $G_q(n-k,n)$. 
In most cases, a universal cycle constructed for $G_q(k,n)$ does not automatically yield one for $G_q(n-k,n)$. 
However, we have observed two exceptional instances where this dual universality does occur:

\begin{itemize}
    \item For $(q,n)=(2,5)$, among all possible orderings of the representatives $c_i$, only two specific exponent sets,
    \[
    (1,4,8,16,3) \quad \text{and} \quad (1,16,8,4,3) \quad (\text{up to cyclic rotation}),
    \]
    produce a cycle that is simultaneously universal for both $G_2(2,5)$ and $G_2(3,5)$.

    \item For $(q,n)=(3,5)$, taking $\alpha$ a root of $1+2x+x^5$, the construction involves two Galois orbits,
    \[
    \{1,3,9,27,81\} \quad \text{and} \quad \{2,6,18,54,162\},
    \]
    where $81$ and $82$ lie in the same $\mathrm{PGL}_2$-orbit.  
    Replacing $81$ by $82$ gives the exponent set
    \[
    \{1,\,54,\,82,\,18,\,2,\,3,\,9,\,162,\,6,\,27\},
    \]
    which yields a cycle that is universal for both $G_3(2,5)$ and $G_3(3,5)$.
\end{itemize}

These examples suggest that dual universality can occasionally be realized when $|G_q(2,5)|=|G_q(3,5)|$, though it appears to be a delicate phenomenon depending on the precise interplay of Galois and $\mathrm{PGL}_2$ orbits.\\

Declaration of generative AI and AI-assisted Technologies in the writing process:  
During the preparation of this work, the authors used ChatGPT (OpenAI) to assist with language editing, terminology refinement, and formatting suggestions. After using this tool, the authors carefully reviewed and edited the content as needed and take full responsibility for the content of the publication.

\bibliographystyle{plain}
\bibliography{reference}
\end{document}